\documentstyle[12pt,leqno]{article}
\begin{document}
\title{Transversal Twistor Spaces of Foliations}
\author{{\normalsize by}\\Izu Vaisman}
\date{}
\maketitle
{\def\thefootnote{*}\footnotetext[1]%
{{\it 1991 Mathematics Subject Classification}
53 C 12, 57 F 30. \newline\indent{\it Key words and phrases}:
Foliations, foliated (projectable)
objects, transversal twistor spaces.}}
\begin{center} \begin{minipage}{12cm}
A{\footnotesize BSTRACT. The transversal twistor space of a foliation
${\cal F}$ of an even codimension
is the bundle ${\cal Z}({\cal F})$
of the complex structures of the fibers of the
transversal bundle of ${\cal F}$.
On ${\cal Z}({\cal F})$ there exists a foliation $\hat{\cal F}$ by covering
spaces of the leaves of ${\cal F}$,
and any Bott connection of ${\cal F}$
produces an ordered pair $({\cal I}_{1}$, ${\cal I}_{2})$
of transversal almost complex
structures of $\hat{\cal F}$.
The existence of a  Bott connection
which yields a structure
${\cal I}_{1}$ that is projectable to the space of leaves is
equivalent to the fact that ${\cal F}$
is a transversally projective foliation. A Bott connection which yields
a  projectable structure ${\cal I}_{2}$ exists iff ${\cal F}$
is a transversally projective foliation which satisfies a supplementary
cohomological condition, and, in this case, ${\cal I}_{1}$ is
projectable as well. ${\cal I}_{2}$ is never integrable.
The essential integrability condition of ${\cal I}_{1}$
is the flatness of the transversal projective structure of
${\cal F}$.
}
\end{minipage} \end{center}
\noindent
The twistor construction \cite{{Pn},{AHS},{OR}} associates (almost)
complex manifolds to differentiable manifolds endowed with a connection,
and, in particular, to Riemannian manifolds with the Levi-Civita connection.
In some cases, this allows to associate objects of Riemannian
(conformal)
geometry to holomorphic objects
of the twistor space. Our aim is to give a similar construction
which will associate a transversally (almost) complex foliation
to a given differentiable foliation. We will define a twistor space which
consists of all the complex structures of the fibers of the transversal bundle
of the foliation. There exists a well defined lift of the original
foliation to the twistor space, the leaves of which are covering spaces of the
original leaves.
Then, a Bott connection of the original foliation
defines two transversal almost
complex structures of the lifted foliation,
and we study the conditions which ensure that
these structures are projectable to the space of leaves, then, the
integrability conditions.

The projectability problem does not exist for twistor
spaces of manifolds.
We show that Bott connections which yield projectable
structures exist iff the foliation has a transversal projective structure which
satisfies a certain cohomological condition.
The general transversally projective foliations are exactly the foliations
which admit a torsionless Bott connection such that a specified one of the two
almost complex structures of the twistor space mentioned above
is projectable.

The integrability conditions are
analogous to those of the classical twistor spaces e.g., \cite{OR}.
In particular, the essential integrability condition of
the projectable almost complex structures transversal to the lift of
a transversally projective foliation is the flatness of the projective
transversal structure
of the foliation. Equivalently, the flatness
condition means that the
foliation has a
Haefliger cocycle which consists of projective transformations expressed
in non homogeneous coordinates.

{\it Acknowledgements}. Part of this work was done during visits of the author
at Istituto di Matematica,
Universit\`a di Roma 1 and Dipartimento di Matematica Universita della
Basilicata, Potenza, Italy
(visit sponsored by the Consiglio Nazionale
delle Ricerche, Italy), Centre de Math\'ematiques,
\'Ecole Polytechnique, Palaiseau, France, and the universities of Jassy and
Bra\c sov (Romania).
The author wishes to express here
his gratitude to these institutions and to his hosts
Paolo Piccinni, Sorin Dragomir, Yvette Kosmann-Schwarzbach,
Fran\c cois Laudenbach, Paul Gauduchon, Radu Miron, Vasile Cruceanu, Vasile
Oproiu, Mihai Anastasiei and Gheorghe Munteanu.
\section{Recall on foliations}
We assume that the reader is familiar
with the differential geometry of foliations, as exposed for instance in
\cite{Mol}, and only recall a few things which we particularly need.
Throughout the paper, {\em differentiable} means
of class $C^{\infty}$, and all the involved objects are of this class.
We are looking at a given {\em foliated manifold} $(M,{\cal F})$ where $M$ is a
manifold of dimension $n=p+s$, and ${\cal F}$ is a foliation
of $M$, of codimension $s$.
We denote by $F$ the tangent bundle $T{\cal F}$, and we
identify the transversal bundle $N{\cal F}=TM/F$ with a fixed complementary
subbundle $E$ (i.e., $TM=E\oplus F$), by sending each vector $Z\in E_{x}$
to the class $[Z]_{F}=\{Z+F_{x}\}\in N_{x}{\cal F}$ $(x\in M)$.
In principle, one is interested in constructions which
depend on $N{\cal F}$, and not
on the choice of $E$ but, the choice of $E$ is
technically helpful since it allows for a natural bigrading of differential forms
and tensor fields. Our convention is that if a bidegree is $(u,v)$
then $u$ is the $E$-degree and $v$ is the (leafwise) $F$-degree.

The most important fact is the existence of a decomposition of the exterior
differential $d$ namely,
$$d=d'_{(1,0)}+d''_{(0,1)}+\partial_{(2,-1)}, \leqno{(1.1)}$$
where the indices denote the bidegrees of the operators.
The operator $d''$ essentially is the exterior differential
along the leaves of ${\cal F}$, and yields the ${\cal F}$-tangential
de Rham cohomology,
which is that of the sheaf
$\Phi(M,{\cal F})$ of germs of functions on $M$ which are
constant along the leaves of ${\cal F}$. On the other hand, a differential
form
of type $(u,0)$ will be called a {\em transverse form}.

It will be our general convention to say that a geometric object on $M$
is {\em ${\cal F}$-foliated
or ${\cal F}$-projectable} (the two terms will be used interchangeably)
if it projects to the local {\em slice spaces} of ${\cal F}$.
For instance,
$f\in C^{\infty}(M)$ is foliated if it is constant along the leaves of
${\cal F}$
i.e., $d''f=0$, a mapping $\varphi:(M_{1},{\cal F}_{1})\rightarrow
(M_{2},{\cal F}_{2})$ is foliated if it sends leaves of ${\cal F}_{1}$
into leaves of ${\cal F}_{2}$. A differential form $\lambda$ on $M$
is foliated if it is the pullback of local forms of the slice spaces
(hence, $\lambda$ also is a transverse form).
A vector field $X$ on $M$ is foliated if it is projectable to the slice spaces
or, equivalently, its flow preserves ${\cal F}$
(see, for instance,
\cite{{Mol},{V2},{V1}} for details).

A principal bundle $P\rightarrow M$ of structure group $G$
is {\em foliated} if it is endowed with a {\em foliated structure},
i.e., a maximal local trivialization atlas
with foliated transition functions. A bundle may have many foliated
structures or none. Any associated bundle of a principal foliated bundle
is again called a {\em foliated bundle}. This definition applies to
vector bundles, in particular. It is an immediate consequence
of the definition that any foliated bundle over $(M,{\cal F})$
has a canonical foliation
$\tilde{\cal F}$, the leaves of which are covering spaces of
the leaves of ${\cal F}$. For the principal bundle $P$, this amounts
to the definition given in \cite{KT}.
Finally, cross sections of a foliated bundle are foliated
if they send leaves of ${\cal F}$ to leaves of $\tilde{\cal F}$.

If $P$ is a foliated principal
bundle over $(M,{\cal F})$, an ${\cal F}$-{\em adapted connection}
on $P$ ({\em Bott connection} \cite{Bt}) is a connection defined
by a ${\cal G}$-valued connection form,
where ${\cal G}$ is the Lie algebra
of $G$, which is of the $\tilde{\cal F}$-type $(1,0)$ (i.e., a transverse
form).
Equivalently, the restricted holonomy of the connection along
paths in the leaves of ${\cal F}$ is trivial.
Furthermore, a connection
which is defined by a foliated connection form is a
{\em projectable connection} (details in \cite{{KT},{Mol}}).

A notion which will be important for us is that of a {\em
projectable (foliated)
distribution}. By definition, a tangent distribution
$D$ of $M$ is projectable if $F\subseteq D$, and, locally, $D$ is
projection-related with distributions of the slice spaces. We have
\proclaim 1.1 Proposition. Let $D$ be a distribution which includes $F$.
Then $D$ is projectable iff it satisfies one of the following
equivalent conditions: a) there exists an open covering of $M$,
$M=\cup_{\alpha}U_{\alpha}$ such that $D/_{U_{\alpha}}$ is spanned by
foliated vector fields $\{X_{k}\}$; b) there exists an open
covering $M=\cup_{\alpha}U_{\alpha}$
such that $D/_{U_{\alpha}}$ has equations
$\xi^{h}_{(\alpha)}=0$, where $\xi^{h}_{(\alpha)}$ are $(1,0)$-forms
which satisfy the condition
$$d''\xi^{h}_{\alpha}=0\hspace{5mm}(mod.\;\xi^{h}_{\alpha}). \leqno{(1.2)}$$
\par \noindent
{\bf Proof.} The equivalence of condition a) with the
definition is obvious.
For condition b), and if $dim\,D=m\geq p$,
it suffices to consider independent forms $\xi^{h}$,
$h=1,...,s+p-m$. If $D$ is projectable, on adapted coordinate
neighborhoods $U_{\alpha}$ of ${\cal F}$ (e.g., \cite{Mol}),
$D$ has equations $\xi^{h}_{(\alpha)}=0$, which are pullbacks
of the equations of the projection of $D$ onto $U_{\alpha}/{\cal F}$.
Then, the forms $\xi^{h}_{(\alpha)}$ are foliated, and (1.2) holds.

Conversely, if $D$ has equations which satisfy (1.2),
any equivalent equations  $\eta^{l}=0$
$(l=1,...,s+p-m)$ where
$\eta^{l}=\sum_{h}A^{l}_{h}\xi^{h}$ also satisfy (1.2). In particular,
we may use a system of the form
$$\eta^{l}=\sigma^{l}+\sum_{u=s+p-m+1}^{s}f^{l}_{u}\sigma^{u},
\leqno{(1.3)}$$
where $(\sigma^{l},\sigma^{u})$ is a foliated basis of $E^{*}$, the dual
bundle of the transversal bundle $E$ of $F$.
Then (1.2) for $\eta^{l}$ means
$$d''\eta^{l}=\sum_{u}d''f^{l}_{u}\wedge\sigma^{u}=
\sum_{h}\eta^{h}\wedge\varphi^{l}_{h}, \leqno{(1.4)}$$
where $\varphi^{l}_{h}$ are some $(0,1)$-forms.
But, the presence of the terms $\sigma^{l}$ in (1.3) make
the last equality (1.4)
impossible unless $\varphi^{l}_{h}=0$. Therefore, $\eta^{l}$ are
foliated forms, and $D$ is foliated. Q.e.d.
\section{Transversal twistor spaces}
In this section, ${\cal F}$ is a foliation of an even
codimension $s=2q$, $q\geq1$, on a manifold $M^{n}$ $(n=p+2q)$.
Similarly to the classical theory
\cite{OR}, we define the {\em transversal twistor
space}
$${\cal Z}({\cal F})=\{J_{x}\in End(N_{x}{\cal F})\;/\;J^{2}=-Id,
\;x\in M\}\leqno{(2.1)}$$
$$\approx\{J_{x}\in End\,E_{x}\;/\;J^{2}=-Id,\,x\in M\},$$
where the equivalence follows from the identification
$N{\cal F}\approx E$ defined in Section 1.
We will study the structure of this space by using a moving frame method
similar to that of \cite{{JR},{V5}}. We express everything in the
$E$-version of (2.1) but, the results may be transposed to the
$N{\cal F}$-version.

Denote by $E^{c}=E\otimes_{\bf R}{\bf C}$ the complexification
of $E$, and consider bases  (frames)
of the fibers $E^{c}_{x}$ $(x\in M)$ which are of
the form $B=(b,\bar b)$, where $b$ is a line matrix of $q$ vectors $(b_{i})
_{i=1}^{q}$, and the bar denotes complex conjugation. (Our convention for
matrix entries is: the upper index is the line index and the lower index is
the column index.) Such bases are obtained by considering $q$-planes
$S\subset E^{c}$ of {\em real index} $0$ (for $S\subset E^{c}$, the real
index is defined as
the dimension of $S\cap\bar S$), and taking bases $b$ of $S$.
Bases $B$ as above will be called {\em non-real bases (frames)}.

We denote by ${\cal B}({\cal F})$ the $Gl(2q,{\bf R})$-principal bundle
of all the non-real frames of the fibers of $E^{c}$, where overall
in this paper
$$Gl(2q,{\bf R})=\{\Phi=\left( \begin{array}{cc}
A&B\\\bar B&\bar A \end{array}\right)\;/\;A,B\in{\cal M}(q,q;{\bf C}),
\;det\, \Phi\neq0\},\leqno{(2.2)}$$
where ${\cal M}(i,j;{\bf C})$
is the space of complex $(i,j)$-matrices.
This obviously is a foliated bundle with the foliated structure defined by
that of
$N{\cal F}$, which is known to be
a foliated vector bundle \cite{Mol}. Therefore,
${\cal B}({\cal F})$  has the lifted foliation
$\tilde{\cal F}$ noticed in Section 1.

Now, we note the existence of the natural projection
$$P:{\cal B}({\cal F})\rightarrow{\cal Z}({\cal F})\leqno{(2.3)}$$
defined by $P(b,\bar b)=J$ where $J$ is the complex structure with the
$(\sqrt{-1})$-eigenspace $span\{b_{i}\}$.
(2.3) is a principal fibration of structure group
$$Gl(q,{\bf C})=\{\left( \begin{array}{cc} A&0\\ 0&\bar A \end{array}\right)
\;/\;A\in{\cal M}(q,q;{\bf C}),\;det \,A\neq0\}. \leqno{(2.4)}$$
The natural projection
$$Q:{\cal Z}({\cal F})\rightarrow M \leqno{(2.5)}$$
is a bundle associated with $\pi:{\cal B}({\cal F})
\rightarrow M$, with the
fiber equal to the homogeneous
space $Gl(2q,{\bf R})/Gl(q,{\bf C})$ of the real dimension $2q^{2}$.
In particular, ${\cal Z}({\cal F})$ is a differentiable manifold of
dimension $p+2q(q+1)$. Moreover, ${\cal Z}({\cal F})$ admits a lift
of the foliation ${\cal F}$ which will be denoted by $\hat{\cal F}$,
and the leaves of $\hat{\cal F}$ are covering spaces of the leaves
of ${\cal F}$.

We discuss geometric structures on ${\cal Z}({\cal F})$ by using some
natural local cotangent bases.

We start with the principal bundle ${\cal B}(M,{\cal F})$ of the tangent
bases of $M$ which are of the form $(B,C)$, where $C$ is a real basis of
$F$ and $B$ is a non-real basis of $E^{c}$ ({\em transversally non-real,
${\cal F}$-adapted bases}). The structure group of
${\cal B}(M,{\cal F})$ is
$$G:=\{\left(\begin{array}{cc} {\cal X}&0\\ 0&{\cal S}\end{array}\right)
\;/\;{\cal X}\in Gl(2q,{\bf R}),\;{\cal S}\in Gl(p,{\bf R})\},\leqno{(2.6)}$$
and there exists a principal fibration
$$K:{\cal B}(M,{\cal F})\rightarrow{\cal B}({\cal F}) \leqno{(2.7)}$$
of group
$$G_{0}:=\{\left( \begin{array}{cc} Id&0\\ 0&{\cal S}\end{array}
\right)\}\subseteq G, \leqno{(2.8)}$$
where $K(B,C)=B$.

On ${\cal B}(M,{\cal F})$, we have the {\em canonical $1$-form} $\theta$
\cite{KN},
where, for $\Xi\in T_{(B,C)}{\cal B}(M,\\
{\cal F})$, $\theta(\Xi)$ is the
column of the components of $\pi_{*}(\Xi)$ with respect to the basis
$(B,C)$, $\pi:{\cal B}(M,{\cal F})\rightarrow M$ again being the natural
projection. We will denote
$$\theta=\left(\begin{array}{c} \beta\\ \bar\beta\\ \gamma
\end{array}\right),\leqno{(2.9)}$$
where $\beta$ is a column $(\beta^{i})$ of complex valued $1$-forms and
$\gamma$ is a column $(\gamma^{u})$ of real valued $1$-forms,
$i=1,...,q;\,u=1,...,p$.
If $g\in G$, the right translation $R_{g}$ on ${\cal B}(M,{\cal F})$acts on
$\theta$ by \cite{KN}
$$R_{g}^{*}\theta=g^{-1}\circ\theta.\leqno{(2.10)}$$
In particular, if $g\in G_{0}$ defined by (2.8), $\beta^{i}$ and
$\bar\beta^{i}$ are invariant by $R_{g}$, and their pullbacks
$s^{*}\beta^{i}$, $s^{*}\bar\beta^{i}$ by a system $\{s\}$ of trivializing
local cross sections of $K$ yield global $1$-forms, again denoted
$\beta^{i},\bar\beta^{i}$ on ${\cal B}({\cal F})$.
Hereafter, it will be our general
convention not to write the pullback mappings, letting the context
indicate the manifolds where the forms are to be considered.

We also notice that
$\tilde\beta:=\left(
\begin{array}{c} \beta\\  \bar\beta \end{array}\right)$
coincides with the transversal canonical $1$-form defined in
\cite{Mol}, and that the right translation $R_{h}$ of ${\cal B}({\cal F})$,
$h\in Gl(2q,{\bf R})$ acts on $\tilde\beta$ by
$$R_{h}^{*}\tilde\beta=h^{-1}\circ\tilde\beta.\leqno{(2.11)}$$
For $\gamma$ of (2.9), if $g\in G_{0}$, then
$$R_{g}^{*}\gamma={\cal S}^{-1}\circ\gamma.\leqno{(2.12)}$$

The forms $\beta^{i},\bar\beta^{i},\gamma^{u}$ (i.e.,
$s^{*}\beta^{i},s^{*}\bar\beta^{i},s^{*}\gamma^{u}$
for local trivializations $s$ of $K$)
will be a part of a local basis of
$T^{*}{\cal B}({\cal F})$.
A complete basis of $T^{*}{\cal B}({\cal F})$
may be obtained by adding an adapted
connection
form $\varpi$ on ${\cal B}({\cal F})$ with values in the Lie algebra
$gl(2q,{\bf R})$ associated with
the representation (2.2). With the
natural symmetric decomposition \cite{OR}
$$gl(2q,{\bf R})=gl(q,{\bf C})\oplus{\cal N},\leqno{(2.13)}$$
and using $(q,q)$-blocks,
we will write
$$\varpi=\left(\begin{array}{cc} \omega&\theta\\
\bar\theta&\bar\omega \end{array} \right)=
\left( \begin{array}{cc} \omega&0\\ 0&\bar\omega \end{array} \right)+
\left( \begin{array}{cc} 0&\theta\\ \bar\theta&0 \end{array} \right).
\leqno{(2.14)}$$
The form of the matrices in the right hand side of (2.14) is characteristic
for the two terms of the decomposition (2.13), and provides the
definition of these two terms.

The reason for taking an adapted connection is that, in this case,
the lifted foliation $\tilde{\cal
F}$ of ${\cal B}({\cal F})$ has the equations
$$\beta=0,\bar\beta=0,\varpi=0,\leqno{(2.15)}$$
and the horizontal distribution of the connection $\varpi$
contains $T\tilde{\cal F}$.

Furthermore, we recall that \cite{KN}
$$R_{h}^{*}\varpi=ad\,h^{-1}\varpi\hspace{5mm}(h\in Gl(2q,{\bf R})).
\leqno{(2.16)}$$
If $h\in Gl(q,{\bf C})$ is given by the matrix (2.4),
(2.16)
reduces to
$$R_{h}^{*}\omega=A^{-1}\omega A,\hspace{2mm} R^{*}_{h}\theta=A^{-1}\theta\bar
A.\leqno{(2.17)}$$
The interpretation of the formulas (2.17) for the principal fibration $P$
of (2.3)
is: (i) $\left( \begin{array}{cc} \omega&0\\ 0&\bar\omega
\end{array}
\right)$ is a connection of $P$, called the {\em associated connection}
of $\varpi$;
(ii) $\left( \begin{array}{cc} 0&\theta \\
\bar\theta&0 \end{array}\right)$ is a horizontal, tensorial $1$-form
\cite{KN} called the
{\em associated tensorial form} of the fibration $P$.

Notice also that, if $h\in Gl(q,{\bf C})$ the relation
(2.11) becomes
$$R^{*}_{h}\beta=A^{-1}\beta,\hspace{2mm} R^{*}_{h}\bar\beta=\bar
A^{-1}\bar\beta.\leqno{(2.18)}$$

We summarize the results in
\proclaim 2.1 Proposition. i). The pullbacks of the entries of the matrices
$(\gamma,\beta,\bar\beta,\theta,\bar\theta)$ of (2.9), (2.14) by
the
(not explicitly written) local cross sections of
the fibration $P$ of (2.3) are
local bases of
the cotangent bundle $T^{*}{\cal Z}({\cal F})$.
ii). The local equations of the lifted foliation $\hat{\cal F}$
of ${\cal F}$ to ${\cal Z}({\cal F})$ are
$$\beta=0,\bar\beta=0,\theta=0,\bar\theta=0,\leqno{(2.19)}$$
and the equations $$\gamma=0 \leqno{(2.20)}$$
yield a well defined complementary
subbundle $\hat E$ of $\hat F=T\hat{\cal F}$ i.e., $T{\cal Z}({\cal F})=
\hat E\oplus\hat F$.
iii). There exists a
decomposition $T{\cal Z}({\cal F})={\cal H}\oplus{\cal V}$,
where ${\cal V}$ is the $2q^{2}$-dimensional
{\em vertical
distribution}, tangent to the fibers of $Q$ (see (2.5)),
and ${\cal H}$ is the $(p+2q)$-dimensional
{\em horizontal distribution}
given by the projection by $P$ of the horizontal distribution of $\varpi$.
The terms of this decomposition have the following equations
$$({\cal H})\;\theta=0,
\bar\theta=0;\;\;({\cal V})\;\beta=0,\bar\beta=0,\gamma=0,
\leqno{(2.21)}$$ and $\hat F\subseteq{\cal H}$.

We must emphasize that the cobases defined in Proposition 2.1 depend on the
choice of the adapted connection $\varpi$.
A change of this connection is of the form \cite{KN}
$$\varpi'=\varpi+\left( \begin{array}{cc}
\xi&\zeta\\ \bar\zeta&\bar\xi \end{array}\right),\leqno{(2.22)}$$
where the additional term is a $gl(2q,{\bf R})$-valued, horizontal,
tensorial $1$-form on ${\cal B}({\cal F})$, which satisfies (2.16).
Hence, in the new cobasis we have the same $\gamma,\beta,\bar\beta$,
but,
$$\theta'=\theta+\zeta,\;\bar\theta'=\bar\theta+\bar\zeta,
\leqno{(2.23)}$$
where $\zeta$ satisfies the second condition (2.17), and vanishes on ${\cal
V}$.
\section{Almost complex structures on ${\cal Z}({\cal F})$}
By Proposition 2.1 i), connection dependent
geometric structures of the twistor space ${\cal Z}({\cal
F})$ may be discovered by looking for geometric objects defined by the cobases
$(\gamma,\beta,\bar\beta,\theta,\bar\theta)$ which are invariant by
transformations of the form (2.12), (2.17), (2.18).
Then, the effect of a change of the connection may be
studied by using (2.23).

Following usual twistor theory, as in \cite{JR} for instance , we notice
the existence of the following
well defined subbundles of the complexification $\hat E^{c}$
of the $\hat{\cal F}$-transversal bundle $\hat E$ of Proposition 2.1
$$({\cal C}_{1})\;\;\gamma=0,\beta=0,\theta=0,\leqno{(3.1)}$$
$$({\cal C}_{2})\;\;\gamma=0,\beta=0,\bar\theta=0.\leqno{(3.2)}$$
Obviously, these are subbundles of real index zero hence, there
exist two well
defined, transversal, almost complex structures of $\hat{\cal F}$, say
${\cal I}_{1},{\cal I}_{2}$, which have ${\cal C}_{1},{\cal C}_{2}$
as their $(-\sqrt{-1})$-eigenspace, respectively
(the sign was chosen to agree with \cite{JR}).
In the $N{\cal F}$-version of (2.1) ${\cal C}_{1},{\cal C}_{2}$ become the
subbundles of $N^{c}\hat{\cal F}:=N\hat{\cal F}
\otimes{\bf C}$ which have the annihilator
$$\beta=0,\theta=0,\leqno{(3.1')}$$
$$\beta=0,\bar\theta=0,\leqno{(3.2')}$$
respectively, in $(N^{c}\hat{\cal F})^{*}\subseteq T^{*}{\cal Z}({\cal F})$,
and the structures
${\cal I}_{1},{\cal I}_{2}$ of $N\hat{\cal F}$
depend on the connection $\varpi$ but, they
do not depend on the choice of $E$.

For a later use, let us also notice that the equations (3.1$'$), (3.2$'$)
define the distributions ${\cal C}'_{1}:=F\oplus{\cal C}_{1}$,
${\cal C}'_{2}:=F\oplus{\cal C}_{2}$, respectively, on ${\cal Z}({\cal F})$.

We will study conditions
which ensure that $\hat{\cal F}$ is a transversally
holomorphic foliation on ${\cal Z}({\cal F})$. A first step, which
is specific for foliations, and is not
encountered in twistor theory on manifolds,
is to find the conditions for ${\cal I}_{1},{\cal I}_{2}$ to be $\hat{\cal
F}$-projectable.

In order to formulate and prove the result
we need some notation and
preparations. We will denote by $\nabla$ the covariant derivative defined on
$E$ by $\varpi$, as well as the equivalent covariant derivative defined on
$N{\cal F}$ by the identification $N{\cal F}\approx E$,
by $R_{\nabla}$ the corresponding curvature
operator,
and by
$$T_{\nabla}(X,Y):=\nabla_{X}(\pi_{E}Y)-\nabla_{Y}(\pi_{E}X)-\pi_{E}[X,Y]
\hspace{5mm}(X,Y\in\Gamma TM),$$
where $\pi_{E}:TM\rightarrow E$ is the projection of the splitting
$TM=F\oplus E$ and $\Gamma$ denotes the space of global cross sections of a
vector bundle,
the torsion \cite{{KN},{Mol}}.
Furthermore,
we may define the {\em Ricci tensor} of an
adapted connection by
$$Ric_{(\nabla,E)}(X,Y)=
tr[Z\mapsto R_{\nabla}(Z,X)\pi_{E}Y]\hspace{3mm}(X,Y\in T_{x}M,
Z\in E_{x}, x\in M),$$
This Ricci tensor depends on the choice of $E$, except if either
$Y\in F_{x}$, hence
$Ric_{(\nabla,E)}(X,Y)=0$, or
$X\in F_{x}$, and then
$$Ric_{(\nabla,E)}(X,Y)=
tr[\pi_{N}Z\mapsto R_{\nabla}(Z,X)\pi_{N}Y]\hspace{3mm}(X,Y,Z\in T_{x}M),$$
where $\pi_{N}:TM\rightarrow{\cal F}$ is the natural projection. In these cases,
we will denote the tensor by $Ric_{\nabla}$.

The {\em structure equations} of the connection $\varpi$ \cite{{KN},{Mol}} are
$$\left( \begin{array}{c} d\beta\\ d\bar\beta \end{array} \right)=
- \left(\begin{array}{cc} \omega&\theta\\ \bar\theta&\bar\omega \end{array}
\right)\wedge\left( \begin{array}{c}\beta\\ \bar\beta \end{array}\right)
+\left(\begin{array}{c}T\\ \bar T\end{array} \right),\leqno{(3.3)}$$
$$\left( \begin{array}{cc} d\omega&d\theta\\ d\bar\theta&d\bar\omega
\end{array}\right)=
-\left( \begin{array}{cc} \omega&\theta\\ \bar\theta&\bar\omega
\end{array}\right)\wedge
\left( \begin{array}{cc} \omega&\theta\\ \bar\theta&\bar\omega
\end{array}\right)+
\left( \begin{array}{cc} \Omega&\Theta\\ \bar\Theta&\bar\Omega
\end{array}\right), \leqno{(3.4)}$$
where the last term of (3.3) consists of the {\em torsion forms},
and the last term
of (3.4) consists of the {\em curvature forms} of the connection.
The forms in (3.3), (3.4) are global forms on
${\cal B}({\cal F})$. Since an adapted connection has no restricted
leafwise holonomy (Section 1),
the curvature forms
vanish if evaluated on two arguments tangent to $\tilde{\cal
F}$. The torsion forms vanish if
at
least one argument is tangent to $\tilde{\cal F}$.

By an exterior differentiation of (3.3) and (3.4), one deduces the
{\em Bianchi identities}. In particular, from the Bianchi identity provided
by (3.3) we get
$$R_{\nabla}(X,Y)Z-R_{\nabla}(Z,Y)X=\nabla_{Y}T_{\nabla}(X,Z)
\hspace{3mm}(Y\in F_{x},X,Z\in E_{x}),\leqno{(3.5)}$$
with the corresponding $N{\cal F}$-version
$$R_{\nabla}(X,Y)\pi_{N}Z-R_{\nabla}(Z,Y)\pi_{N}X\leqno{(3.5')}$$
$$=\nabla_{Y}T_{\nabla}(X,Z)
\hspace{3mm}(Y\in F_{x},X,Z\in T_{x}M).$$
We will denote by $i(X)\nabla_{Y}T_{\nabla}$ the endomorphism
$Z\mapsto \nabla_{Y}T_{\nabla}(X,Z)$ of $E$,
and notice the following consequence of (3.5$'$)
$$Ric_{\nabla}(Y,X)=tr[R_{\nabla}(X,Y)-i(X)\nabla_{Y}T_{\nabla}]
\hspace{3mm}(Y\in F_{x},X\in T_{x}M).\leqno{(3.5'')}$$

Now, we can state the following results

\proclaim 3.1 Theorem. 1).
The $\hat{\cal F}$-transversal
almost complex structure
${\cal I}_{1}$ is projectable iff  either i) $q\geq2$, and
$\forall x\in M$, $\forall Y\in F_{x}$, $\forall X,Z\in T_{x}M$, one has
$$R_{\nabla}(X,Y)\pi_{N}Z=\lambda(X,Y)\pi_{N}Z+\mu(Z,Y)\pi_{N}X,\leqno{(3.6)_{1}}$$
where, necessarily, the coefficients $\lambda,\mu$ are given by
$$\lambda(X,Y)=\frac{1}{4q^{2}-1}(2q\,tr\,R_{\nabla}(X,Y)-Ric_{\nabla}(Y,X))
\leqno{(3.6')_{1}}$$
$$=\frac{1}{2q+1}tr\,R_{\nabla}(X,Y)+
\frac{1}{4q^{2}-1}tr\,i(X)\nabla_{Y}T_{\nabla},$$
$$\mu(X,Y)=\frac{1}{1-4q^{2}}(tr\,R_{\nabla}(X,Y)-2q\,Ric_{\nabla}(Y,X))
\leqno{(3.6'')_{1}}$$
$$=-\frac{1}{2q+1}tr\,R_{\nabla}(X,Y)-
\frac{2q}{4q^{2}-1}tr\,i(X)\nabla_{Y}T_{\nabla}$$
or ii) $q=1$, and
$$R_{\nabla}(X,Y)\pi_{N}X=\nu(X,Y)\pi_{N}X,\leqno{(3.6)_{2}}$$
where $X,Y$ are as above and necessarily
$$\nu(X,Y)
=\frac{1}{3}[tr\,R_{\nabla}(X,Y)+Ric_{\nabla}(Y,X)]\leqno{(3.6')_{2}}$$
$$=\frac{2}{3}tr\,R_{\nabla}(X,Y)-\frac{1}{3}tr\,
i(X)\nabla_{Y}T_{\nabla}.$$
2). For all $q$, the structure
${\cal I}_{2}$ is projectable
iff
$$R_{\nabla}(X,Y)=\alpha(X,Y)Id,\leqno{(3.6)_{3}}$$
where $Id$ is the identity mapping of $N{\cal F}$ and, necessarily
$$\alpha(X,Y)=\frac{1}{2q}tr\,R_{\nabla}(X,Y),\leqno{(3.6')_{3}}$$
and in this case the structure
${\cal I}_{1}$ must also be projectable.
\par
\noindent{\bf Proof.} While the results are formulated for $N{\cal F}$,
for the proof we choose a transversal bundle $E$, and take $X,Z\in E_{x}$,
forgetting about the $\pi_{N}$. Generally,
if $(V,{\cal S})$ is a foliated manifold,
projectability of an almost complex structure $J$ of
$N{\cal S}$ onto the local slice spaces,
is equivalent to the fact that for any
foliated
cross section $s$ of $N{\cal S}\approx E$
($E$ is transversal to $T{\cal S}$), $Js$ is also foliated, as well as
to the fact
that, if ${\cal C}\subseteq E^{c}$ is
the $(-\sqrt{-1})$-eigendistribution
of $J$, then ${\cal C}':=T{\cal S}\oplus{\cal C}$ is projectable.
In our case, we will discuss this property for the distributions
${\cal C}'_{1},{\cal C}'_{2}$ defined by (3.1$'$), (3.2$'$).

Since the curvature and torsion forms are horizontal on
the
principal bundle ${\cal B}({\cal F})$, we may write
the entries of the torsion and curvature matrices as follows
$$T^{a}=\frac{1}{2}T^{a}_{bc}\beta^{b}\wedge\beta^{c}
+T^{a}_{b\bar c}\beta^{b}\wedge\beta^{\bar c}
+\frac{1}{2}T^{a}_{\bar b\bar c}\beta^{\bar b}\wedge\beta^{\bar c},
\leqno{(3.7)}$$
$$\Theta^{a}_{\bar b}=\frac{1}{2}R^{a}_{\bar bcd}\beta^{c}\wedge\beta^{d}
+R^{a}_{\bar bc\bar d}\beta^{c}\wedge\beta^{\bar d}
+\frac{1}{2}R^{a}_{\bar b\bar c\bar d}\beta^{\bar c}\wedge\beta^{\bar d}
\leqno{(3.8)}$$
$$+
R^{a}_{\bar bcu}\beta^{c}\wedge\gamma^{u}+
R^{a}_{\bar b\bar cu}\beta^{\bar c}\wedge\gamma^{u},$$ etc.
In formulas (3.7), (3.8), we use the Einstein summation convention,
$a,b,c,d=1,...,q$, and $u=1,...,p$. The coefficients, except those
in
the last two terms of (3.8), are global functions on
${\cal B}({\cal F})$. The indices in the left hand side of (3.8) are in
agreement with the position of the block $\Theta$ in the curvature matrix.

If we pull back equations (3.3), (3.4), (3.7), (3.8) by local cross
sections $s$
of $P$
of (2.3), we get the same equations on
${\cal Z}({\cal F})$,
locally. The terms of (3.7), (3.8) have corresponding
bidegrees with respect to the decomposition
$T{\cal Z}({\cal F})=\hat E\oplus\hat F$, and we see that,  if $d''$ is the
$\hat{\cal F}$-leafwise differential corresponding to this
decomposition, then
$$d''\beta=0,\;d''\theta=(R^{a}_{\bar bcu}\beta^{c}+R^{a}_{\bar b\bar
cu}\bar\beta^{c})\wedge\gamma^{u}.\leqno{(3.9)}$$

Accordingly, Proposition 1.1 tells us that ${\cal C}'_{1}$ is projectable
iff
$$R^{a}_{\bar b\bar
cu}=0,\leqno{(3.10)}$$
and ${\cal C}'_{2}$ is projectable iff
$$R^{a}_{\bar b
cu}=0.\leqno{(3.11)}$$
Of course, (3.10), (3.11) must hold everywhere on ${\cal Z}({\cal F})$,
therefore, everywhere on ${\cal B}(M,{\cal F})$
too (notation of Section 2).

For any given frame $(B,C)\in {\cal B}_{x}(M,{\cal F})$
($x\in M$), we have the following interpretation of the matrix $\Theta$ as
an $(End\,E_{x})$-valued $2$-form
$$\Theta=pr_{S}\circ R_{\nabla}\circ pr_{\bar S},\leqno{(3.12)}$$
where $S=span\{b\}$, $\bar S=span\{\bar b\}$ $(B=(b,\bar b))$.
Accordingly (3.10) means that, $\forall x\in M$, $\forall Y\in F_{x}$, $\forall
S^{q}
\subseteq E^{c}_{x}$ of real index zero, $\forall V,W\in S$,
$$R_{\nabla}(\bar V,Y)\bar W\in \bar S.\leqno{(3.13)}$$
Indeed, (3.13) implies (3.10) since the latter is (3.13) with arguments
belonging to a tangent basis. On the other hand, (3.10) implies (3.13)
since,
for any given $Y,V,W,S$
there are frames $(B,C)$ such that $S=span\{b\}$,
$V=\alpha b_{c}$, $W=\beta_{1} b_{b}
+\beta_{2}b_{c}$, $Y=\gamma c_{u}$ $(\alpha,\beta_{1},\beta_{2}\in{\bf C},
\gamma\in{\bf R})$.
Furthermore,
for any given $V,W$ with $span\{V,W\}$
of real index zero, there are
subspaces
$S_{1},S_{2}$ of real index zero
such that $S_{1}\cap S_{2}=
span\{V,W\}$.
Hence, (3.13) is equivalent to the fact that,
for some coefficients $\lambda,\mu$,
$$R_{\nabla}(\bar V,Y)\bar W=\lambda(\bar V,Y)\bar W+\mu(\bar W,Y)\bar V,
\leqno{(3.14)}$$ holds
whenever $span\{V,W\}$ has real index zero.

Now, notice that $span\{V,W\}$ is of real index zero iff
either $V=0,W=0$, or
$V,W$ are both proportional to some $X+\sqrt{-1}Z$ where $X,Z\in E$
are ${\bf R}$-independent, or $V,W$ are ${\bf C}$-independent,
and their real and imaginary parts are four ${\bf R}$-independent
vectors. (Use the fact that the real index is zero iff
the real dimension of
$span\{V,\bar V,W,\bar W\}$ is $4$ in the first case, and $8$ in the
second case.) Therefore, if $span\{V,W\}$ is of real index zero, the same holds
if one or both vectors $V,W$  are replaced by their
complex conjugates.

Accordingly, if $q\geq2$, and if (3.14) holds for  $V,W$, (3.14) also
holds for $(\bar V,W)$, $(V,\bar W)$, $(\bar V,\bar W)$, therefore,
(3.14) holds iff the condition (3.6)$_{1}$ holds.

In the case $q=1$, (3.14) reduces to
$$R_{\nabla}(V,Y)V=\nu(V,Y)V\leqno{(3.15)}$$
for any $V=X+\sqrt{-1}Z$ with independent $X,Z\in E_{x}$, and for some
coefficient $\nu$.
Separating the imaginary and real parts, (3.15) becomes
$$R_{\nabla}(X,Y)Z+R_{\nabla}(Z,Y)X=\nu(X,Y)Z+\nu(Z,Y)X,
\leqno{(3.16)}$$
$$R_{\nabla}(X,Y)X-R_{\nabla}(Z,Y)Z=\nu(X,Y)X-\nu(Z,Y)Z,
\leqno{(3.17)}$$
and (3.17) follows from (3.16) by taking $X=Z$.
Finally, (3.16) is just the polarization of the $E$-version of
(3.6)$_{2}$ with respect to $X$.

The similar analysis of (3.11) requires the replacement of $\bar V$ by $V$
in (3.13), and using the same argument as for
(3.14), the projectability
condition of ${\cal C}_{2}$ will be
$$R_{\nabla}(V,Y)\bar W=\alpha(V,Y)\bar W,
\leqno{(3.18)}$$ where the real index of $span\{V,W\}$ is zero, and
$\alpha$ is a corresponding coefficient.
No term in $\bar V$ appears in the right hand side
since $\bar V$ does not appear in the left hand side of (3.18).

Now, for $q\geq2$, again, we may replace the pair $(V,W)$ by any
of the pairs $(\bar V,W)$, $(V,\bar W)$, $(\bar V,\bar W)$,
and we see that (3.18) is equivalent to
$$R_{\nabla}(X,Y)Z=\alpha(X,Y)Z, \leqno{(3.19)}$$
where $X,Z\in E_{x},Y\in F_{x},x\in M$.

If $q=1$, we have to use (3.18) for $W=V=X+\sqrt{-1}Z$, with independent
vectors $X,Z$, and the result is equivalent to the conditions
$$R_{\nabla}(X,Y)Z-R_{\nabla}(Z,Y)X=\alpha(X,Y)Z-\alpha(Z,Y)X,
\leqno{(3.20)}$$
$$R_{\nabla}(X,Y)X+R_{\nabla}(Z,Y)Z=\alpha(X,Y)X+\alpha(Z,Y)Z.
\leqno{(3.21)}$$
But, (3.21) holds iff $R_{\nabla}(X,Y)X=\alpha(X,Y)X$,
which is equivalent to (3.16) by polarization, and the pair of
conditions (3.16), (3.20) is equivalent to
(3.19). Thus, the only projectability condition is (3.19) again,
and we have proven (3.6)$_{3}$.
Since  (3.19) is (3.6)$_{1}$ for $\mu=0$, we see that the
projectability
of ${\cal C}_{2}$ implies the projectability of ${\cal C}_{1}$, if $q\geq2$.
The same is true for $q=1$ since (3.19) implies (3.6)$_{2}$.

Hence, the condition of simultaneous projectability of ${\cal C}_{1},
{\cal C}_{2}$ is condition (3.6)$_{3}$
of the projectability of ${\cal C}_{2}$.

Now, to, end the proof of Theorem 3.1, it
remains to justify the expressions
of the coefficients $\lambda,\mu,\nu,\alpha$ as traces.
Using the following hints,
these are easy consequences of the initial formulas (3.6).
For (3.6$'$)$_{1}$, (3.6$''$)$_{1}$,
compute the traces of the operators which send $Z$, respectively $X$, to
the vector defined by the $E$-version of
(3.6)$_{1}$ then, use (3.5), (3.5$''$).
For (3.6$'$)$_{2}$, polarize the $E$-version of
(3.6)$_{2}$ with respect to $X$ before computing traces.

Notice that if the splitting $TM=E\oplus F$ is chosen
$\lambda,\mu,\nu,\alpha$ may be seen
as forms of type $(1,1)$ on $(M,{\cal F})$. Otherwise, we may see them
as cross sections of $Hom(TM\otimes F,{\bf R})$ which vanish if the
first argument is in $F$.
Q.e.d.

Clearly, if $\nabla$ is an
${\cal F}$-projectable connection, the almost complex structures
${\cal J}_{1},{\cal I}_{2}$ will be $\hat{\cal F}$-projectable.
An interpretation of the simultaneous
projectability condition (3.6)$_{3}$ is given
by
\proclaim 3.2 Proposition. The almost complex structure ${\cal I}_{2}$
(hence,
${\cal I}_{1}$ too), is projectable iff there exists an open
covering $M=\cup U_{\alpha}$ endowed with projectable connections
$\nabla^{\alpha}$ of $N{\cal F}/_{U_{\alpha}}$ such that
$\nabla/_{U_{\alpha}}$ and $\nabla^{\alpha}$ are related by a
semisymmetric transformation. \par
\noindent{\bf Proof.} The structure equations (3.4) show that, locally,
$$R_{\nabla}(X,Y)=d''\varpi_{M}(X,Y),\leqno{(3.22)}$$
where $\varpi_{M}$ is the pullback of $\varpi$ to $M$ by local cross sections
of ${\cal B}({\cal F})$. Hence, locally, (3.6)$_{3}$ may be seen
as
$$d''\varpi_{M}=\alpha\otimes Id=d''\sigma\otimes
Id\hspace{5mm}(\alpha=d''\sigma)\leqno{(3.23)}$$
(the first equality (3.23) implies $d''\alpha=0$; then,
see \cite{V1} for the $d''$-Poincar\'e lemma.) Accordingly, the
local connections defined on $N{\cal F}$ by
$$\theta=\varpi_{M}-\sigma\otimes Id\leqno{(3.24)}$$
are projectable. On the other hand, (3.24) exactly is what the older books
on differential geometry used to call a {\em semisymmetric transformation}
(the notion was used when one of the two connections was
torsionless, which we do not ask here).
Equivalently, it is easy to understand that two connections on $N{\cal F}$
are semisymmetrically related iff they define the same horizontal
distribution ${\cal H}$ (see (2.21)) on the twistor space ${\cal Z}({\cal
F})$.

The converse follows from
\proclaim 3.3 Proposition. Two adapted connections $\varpi,\varpi'$ of
$N{\cal F}$ define the same transverse almost complex structures
${\cal I}_{1},{\cal I}_{2}$ of $\hat{\cal F}$ iff they are related by a
semisymmetric transformation. \par
\noindent{\bf Proof.} The local equations of the distributions (3.1), (3.2)
of the two connections contain the same forms $\gamma,\beta$,
and forms $\theta,\theta'$ which are related by (2.23) with the horizontal
term $\zeta$. Since the connections are adapted,
$\zeta$ has no term in $\gamma$
(notation of Proposition 2.1), and if we want ${\cal C}_{1}={\cal C}'_{1}$,
${\cal C}_{2}={\cal C}'_{2}$, $\zeta$ also cannot have terms
in $\beta,\bar\beta$. Hence, the almost complex structures defined by the
two connections coincide iff $\theta'=\theta$.
With the same notation as in (3.12), this condition has the
equivalent form
$$pr_{S}\circ(\varpi'_{M}-\varpi)_{M}(X)\circ pr_{\bar S}=0,
\hspace{3mm} X\in TM,\leqno{(3.25)}$$
and arguments similar to those in the proof of Theorem 3.1 show
that (3.25) is equivalent to
$$\varpi'=\varpi+\tau\otimes Id,\leqno{(3.26)}$$
where $\tau$ is a $1$-form which vanishes on the leaves of ${\cal F}$.
Q.e.d.

It is also worthwhile noticing
\proclaim 3.4 Proposition. A torsionless adapted connection $\varpi$ defines a
projectable almost complex structure ${\cal I}_{2}$ (therefore, ${\cal I}_{1}$
too) iff $\varpi$ is a projectable connection.\par
\noindent{\bf Proof.} By inserting
(3.6$'$)$_{3}$ in (3.5), and taking the trace, it
follows that $\alpha=0$. Since (3.6)$_{3}$
is locally equivalent to (3.23), $\varpi$ must be a foliated form. Q.e.d.

As a matter of fact, the result is true under the weaker hypothesis
$\nabla_{Y}T_{\nabla}=0$, $\forall Y\in F$.

On the other hand,
the projectability of ${\cal I}_{1}$ alone does not imply the
projectability of $\varpi$. Indeed, if
(3.6)$_{1}$ is inserted
in (3.5) we only get $\lambda=\mu$  i.e.,
if $q\geq2$, a torsionless adapted connection (or one with
$\nabla_{Y}T_{\nabla}=0$) defines a projectable structure ${\cal I}_{1}$
iff its curvature operator is
$$R_{\nabla}(X,Y)Z=\lambda(X,Y)Z+\lambda(Z,Y)X,\hspace{3mm}Y\in F,\;X,Z\in E,
\leqno{(3.27)}$$ where
$$\lambda(X,Y)=\frac{1}{2q+1}tr\,R_{\nabla}(X,Y).\leqno{(3.28)}$$
By polarizing the $E$-version of
(3.6)$_{2}$ with respect to $X$, and using (3.5), we see that the
same is true for $q=1$.

But, if we add the hypothesis that $\nabla$ induces a projectable
connection in $\wedge^{2q}N{\cal F}$, which is equivalent to
$tr\,R_{\nabla}(X,Y)=0$ for $Y\in F,X\in E$, we see that $\nabla$ itself
must be a projectable connection.

We end this section by
\proclaim 3.5 Proposition. On ${\cal Z}({\cal F})$, the horizontal
distribution ${\cal H}$ defined by (2.21) is integrable iff $\varpi$ is
locally semisymmetric flat. \par
\noindent{\bf Proof.} The integrability of ${\cal H}$ means
$d\theta=d\bar\theta=0$ $(mod.\,\theta,\bar\theta)$, and (3.4) tells us
that this happens iff $\Theta=0$. By treating this last condition as we did
with (3.12), the integrability condition becomes
$$R_{\nabla}(X,Y)=\lambda(X,Y)Id,\hspace{5mm}X,Y\in\Gamma
TM,\leqno{(3.29)}$$
where $\lambda$ is a $2$-form on $M$ which is closed because of the Bianchi
identity. Thus, locally, $\lambda=d\varphi$ for a $1$-form $\varphi$, and
(3.29) holds iff the local connections $\varpi_{M}-\varphi\otimes Id$ are
flat. Q.e.d.
\proclaim 3.6 Corollary. The expression
$$\Xi:=\sqrt{-1}tr(\theta\wedge\bar\theta)\leqno{(3.30)}$$
pulls back to a well defined, real, global $2$-form of rank $2q^{2}$ on
${\cal Z}({\cal F})$, which is closed iff $\varpi$ is a locally
semisymmetric flat connection.\par
\noindent{\bf Proof.} That $\Xi$ is well defined follows from (2.17). Then,
$$d\Xi=\sqrt{-1}(\Theta\wedge\bar\theta-\theta\wedge\bar\Theta).
\leqno{(3.31)}$$
Thus, $\Theta=0$ implies $d\Xi=0$. The converse holds since $ker\,\Xi={\cal
H}$, and the kernel of a closed $2$-form is involutive. Q.e.d.
\section{Existence of suitable connections}
Let $(M,{\cal F})$ be a foliated manifold as in the previous sections. If an
adapted connection $\varpi$ of $N{\cal F}$
defines a foliated almost complex structure
${\cal I}_{1}$ on the transversal twistor space
${\cal Z}{\cal F}$
we will say that $\varpi$ is a {\em twistor-suitable
connection}. If $\varpi$ is such that ${\cal I}_{2}$
(hence, ${\cal I}_{1}$ too) is projectable, we will say that
$\varpi$ is a {\em strongly twistor-suitable connection}.
In this section, we study the existence of suitable connections.
\proclaim 4.1 Proposition. Let ${\cal F}$ be an even-codimensional foliation
of a manifold $M$. Then, $N{\cal F}$ has a
strongly twistor-suitable connection iff
the Atiyah class $\alpha({\cal F})\in H^{1}(M,\underline{\wedge^{1,0}M
\otimes End\,N{\cal F}})$ is in the image of
$H^{1}(M,\underline{\wedge^{1,0}M})$ by the mapping $\otimes Id.$
\par
\noindent{\bf Proof.} We recall that the Atiyah class is the obstruction to
the existence of a projectable connection on $N{\cal F}$, and send to
\cite{{Mol2},{Mol}} for details. Underlining in the formulas of the
proposition means that we take the corresponding sheaf of
{\em foliated germs}. $\wedge^{1,0}M$ is the space of differential forms of
the ${\cal F}$-type $(1,0)$.

A representative differential form of
$\alpha({\cal F})$ is given by the part of type $(1,1)$ of the curvature of
any adapted connection $\sigma$ on $N{\cal F}$
\cite{Mol2}. Thus, if the required
connection $\varpi$ exists, and we take $\sigma=\varpi$,
the projectability condition
(3.6)$_{3}$ (see also (3.23)) implies the condition of Proposition 4.1.

Conversely, and with the notation of formula (3.23), the condition of
Proposition 4.1 means
$$d''\sigma_{U_{\alpha}}=\alpha\otimes Id+d''\eta,\leqno{(4.1)}$$
where $M=\cup_{\alpha}U_{\alpha}$
($U_{\alpha}$ open in $M$), $\alpha$ is a global
$d''$-closed form of
${\cal F}$-type (1,1), and $\eta$ is a global
$(End\,N{\cal F})$-valued form of
${\cal F}$-type (1,0). Then, $\varpi=\sigma-\eta$ is a new adapted
connection which defines projectable structures ${\cal I}_{1},{\cal I}_{2}$
since it satisfies condition (3.6)$_{3}$. Q.e.d.

For a more geometric answer, we need the following preparations which apply
some classical notions to adapted connections of $N{\cal F}$.

Let $\nabla^{1},\nabla^{2}$ be the covariant derivatives of two adapted
connections of $N{\cal F}$. Then, there exists a {\em difference tensor}
$$S_{x}(X,Y):=\nabla^{2}_{X}(\pi_{N}\tilde Y)-\nabla^{1}_{X}(\pi_{N}\tilde Y)
\hspace{3mm}(X,Y\in T_{x}M,\,x\in M),\leqno{(4.2)}$$
where $\tilde Y$
is a vector field such that $\tilde Y(x)=Y$.
In particular, by using a foliated field $\tilde Y$, we see that
$S$ vanishes if any
of its arguments is in $F$.
Formula (3.24) tells that $\nabla^{1},\nabla^{2}$ are
semisymmetrically related iff
$$S=\sigma\otimes \pi_{N}\leqno{(4.3)}$$
for some scalar $(1,0)$-form $\sigma$ on $(M,{\cal F})$.

On the other hand, we need the
{\em transposed connection} of a connection $\nabla$, defined by
$$^{t}\nabla_{X}(\pi_{N}Y)=\nabla_{Y}(\pi_{N}\tilde X)+\pi_{N}[\tilde X,Y],
\leqno{(4.4)}$$
where $\tilde X(x)=X\in T_{x}M$, $x\in M$, and $\tilde X$ is foliated.
Using $^{t}\nabla$, we get the torsionless connection
$$^{s}\nabla_{X}(\pi_{N}Y)=\frac{1}{2}(\nabla_{X}(\pi_{N}Y)+
^{t}\nabla_{X}(\pi_{N}Y)),\leqno{(4.5)}$$
called the {\em symmetric part of} $\nabla$.

Furthermore, two torsionless adapted connections $\nabla^{1},\nabla^{2}$
of $N{\cal F}$ are said to be {\em projectively related} if their difference
tensor
is of the form
$$S=\phi\odot\pi_{N},\leqno{(4.6)}$$
where $\phi$ is a $(1,0)$-form on $M$ and
$\odot$ denotes the symmetric tensor product. The name comes from the fact
that, on manifolds (i.e., ${\cal F}$ is the foliation by points), two such
connections have the same selfparallel lines (e.g., \cite{Kob}).
Accordingly, a foliation ${\cal F}$ is called {\em transversally projective}
\cite{{NS},{V4}}
if it is endowed with the following data: i) an open covering
$M=\cup_{\alpha}U_{\alpha}$ by ${\cal F}$-adapted coordinate neighborhoods,\\
ii) a family of foliated connections $\nabla^{\alpha}$ on $N{\cal F}/_{U_{\alpha}}$
which are projectively related over $U_{\alpha}\cap U_{\beta}$ i.e.,
$$S_{\alpha\beta}:=\nabla^{\beta}-\nabla^{\alpha}=\phi_{\alpha\beta}
\odot\pi_{N}\hspace{3mm}(\phi_{\alpha\beta}\in\wedge^{1}_{fol}(U_{\alpha}\cap
U_{\beta})),\leqno{(4.7)}$$
where the index $fol$ means foliated forms.
If maximal, this data system defines a {\em transversal projective structure}
of ${\cal F}$. In the framework of $\check {\rm C}$ech cohomology,
$\{\phi_{\alpha\beta}\}$ defines a cohomology class $[\phi]\in H^{1}(M,
\underline{\wedge^{1,0}M})$ which we call the {\em complementary class}
of the projective structure. If we put
$$\phi_{\alpha\beta}=\psi_{\alpha}-\psi_{\beta},\leqno{(4.8)}$$
where $\psi_{\alpha},\psi_{\beta}$ are transverse $1$-forms on $U_{\alpha},
U_{\beta}$, we get the global
adapted connection
$$\nabla=\nabla^{\alpha}+\psi_{\alpha}\odot\pi_{N},\leqno{(4.9)}$$
which is projectively related to the local connections
$\nabla^{\alpha}$ but, of
course, it may not be foliated. A global foliated connection (4.9) can be
obtained iff $[\phi]=0$.

Now, assume that $N{\cal F}$ has a
strongly suitable connection $\varpi$, and let
$\{\nabla^{\alpha}\}$ be the associated system of local, foliated connections,
semisymmetrically related to $\varpi$, defined by Proposition 3.2. Then, any two
connections of this system are semisymmetrically related, say
$$\nabla^{\beta}=\nabla^{\alpha}+2\phi_{\alpha\beta}\otimes
\pi_{N},\leqno{(4.10)}$$
and it follows that the symmetric connections $^{s}\nabla^{\alpha}$ satisfy
(4.7), and define a transversal projective structure of ${\cal F}$.

Conversely, assume that ${\cal F}$ is a {\em transversally projective foliation},
with the structure defined by the data of (4.7). Then, the
connections
$$\tilde\nabla^{\alpha}=\nabla^{\alpha}+\frac{1}{2}\tau_{\alpha},
\leqno{(4.11)}$$
where $\nabla_{\alpha}$ are those of (4.7),
have the torsion tensor $\tau_{\alpha}$, and for these
connections one has the transition relations
$$\tilde\nabla^{\beta}=\tilde\nabla^{\alpha}+\phi_{\alpha\beta}\odot
\pi_{N}+\xi_{\alpha\beta},
\hspace{5mm}\xi_{\alpha\beta}:=\frac{1}{2}(\tau_{\alpha}-\tau_{\beta}).
\leqno{(4.12)}$$

We would like to be able to choose $\tau$ such that (4.12)
would be a semisymmetric transformation of foliated transversal
connections.
This happens iff
the local tensors $\tau_{\alpha}$ are foliated, and
$$\xi_{\alpha\beta}=\pi_{N}\wedge\phi_{\alpha\beta}.\leqno{(4.13)}$$
The right hand side of the equality (4.13) defines a
$\check{\rm C}$ech cohomology class
$[\phi]_{N}\in
H^{1}(M,\underline{\wedge^{2,0}M\otimes
N{\cal F}})$, which we call the {\em normal complementary class} of the
projective structure. Then, clearly, we can get the required
situation of foliated, semisymmetric relations iff
$[\phi]_{N}=0$. Therefore, we have proven
\proclaim 4.2 Theorem. Let ${\cal F}$ be an even-codimensional foliation of a
manifold $M$. Then, $N{\cal F}$ has an adapted connection which defines
projectable, almost complex structures ${\cal I}_{1},{\cal I}_{2}$ of $N\hat{\cal F}$
on ${\cal Z}({\cal F})$ iff ${\cal F}$ has a transversal, projectable, projective
structure, with a vanishing normal complementary class.

The existence condition of (not strongly) suitable connections is weaker.
Namely, we have
\proclaim 4.3 Theorem. If there exists an adapted connection of $N{\cal F}$
which produces a projectable almost complex structure ${\cal I}_{1}$ of
$N\hat{\cal F}$, there also exists a torsionless connection with the same
property. A torsionless twistor-suitable connection exists iff ${\cal F}$
is a transversally projective foliation. \par
\noindent{\bf Proof.} The $E$-version of the
projectability condition (3.6)$_{1}$ holds
pointwisely iff it holds for any vector field $Y\in\Gamma F$, and any
projectable vector fields $X,Z\in\Gamma E$. Then, if we express $R_{\nabla}$
by covariant derivatives and use Bott's condition
$$\nabla_{Y}X=\pi_{E}[Y,X],\hspace{5mm}Y\in \Gamma F,\;X\in\Gamma E,$$
(see the $N{\cal F}$-version of this condition in \cite{Bt}),
we get the following equivalent form of (3.6)$_{1}$
$$\pi_{E}[Y,\nabla_{X}Z]=-\lambda(X,Y)Z-\mu(Z,Y)X,\leqno{(4.14)}$$
for all $Y\in\Gamma F$, and for all projectable $X,Z\in\Gamma E$.

Since from the Bott condition and the projectability of $X,Z$ we get
$$\pi_{E}[Y,\pi_{E}[X,Z]]=\nabla_{Y}\pi_{E}[X,Z]=0,$$
the $E$-version of (4.5) yields
$$\pi_{N}[Y,^{s}\nabla_{X}Z]=\frac{1}{2}\pi_{E}([Y,\nabla_{X}Z]
+[Y,\nabla_{Z}X]),$$
and (4.14) implies that the symmetric part of $\nabla$ satisfies (3.27)
with
$(\lambda+\mu)/2$ in the role of $\lambda$.

This proves the first conclusion for $q\geq2$. The proof for $q=1$
is similar but, starting with (3.6)$_{2}$ which is equivalent to
$$\pi_{E}[Y,\nabla_{X}X]=-\nu(X,Y)X,\leqno{(4.15)}$$
for all $Y\in\Gamma F$, and for all foliated $X\in\Gamma E$.

For the second assertion of Theorem 4.3, if ${\cal F}$ is transversally
projective, it has the connection $\nabla$ of (4.9) where
$\nabla^{\alpha}$ are foliated connections. From the $E$-version of
(4.9), it follows that $\nabla$
satisfies (4.14) and (4.15) for $\lambda=\mu=d''\psi_{\alpha}$.
$\lambda$ is a global form since the differences $\psi_{\alpha}-\psi_{\beta}$
are foliated $1$-forms.

Conversely, if $\nabla$ is a twistor-suitable, torsionless connection,
we have (3.27), (3.28) where $d''\lambda=0$ hence, $\lambda=d''\sigma$,
$\sigma$ being a local form of ${\cal F}$-type $(1,0)$.
Then, we may locally write (3.27) as
$$\nabla_{Y}\nabla_{X}Z=(Y\sigma(X))Z+(Y\sigma(Z))X,\leqno{(4.16)}$$
where $Y\in\Gamma F$ and
$X,Z$ are foliated cross sections of $E$. From (4.16), it follows that the
local projectively related connection
$$\nabla'=\nabla-\sigma\odot\pi_{N}$$ is foliated. Q.e.d.

Notice that any connection obtained from $\nabla$ of (4.9) by a projective
transformation also yields a projectable structure ${\cal I}_{1}$. But, the
structures ${\cal I}_{1}$ of the different connections are different.
Indeed, the form $\phi$ of a projective difference (4.6) of connections is
real hence, $\theta$ of (3.1) changes by an additional term which contains
both $\beta$ and $\bar\beta$.

Furthermore, if ${\cal F}$ is a transversally projective foliation with a
vanishing normal complementary class, the connections which it generates by
the procedures of Theorems 4.2, 4.3 may not be semisymmetrically related
hence, produce different structures ${\cal I}_{1}$, ${\cal I}_{2}$.

\proclaim 4.4 Remark. In the same way as for Proposition 4.1, we can see that
a foliation ${\cal F}$ has a torsionless twistor-suitable connection
(equivalently, it has a projectable transversal projective structure) iff the
Atiyah class of ${\cal F}$ is in the image of the mapping
$$\odot Id:H^{1}(M,\underline{\wedge^{1,0}M})\rightarrow
H^{1}(M,\underline{\wedge^{1,0}M\otimes End\,N{\cal F}}).$$

A way of looking for examples could be the use of a Riemannian metric $g$
on $(M,{\cal F})$. Then, we will take $E$ orthogonal to $F$ and there exists
a unique adapted torsionless connection $\nabla$ such that the
$\nabla$-parallelism along paths tangent to $E$ preserves the metric
$g/_{E}$. Namely, $\nabla_{Y}X=\pi_{E}[Y,X]$ if $Y\in\Gamma F,X\in\Gamma
E$, and $\nabla_{Y}X$ is determined like the Levi-Civita connection
\cite{KN} if $Y,X\in\Gamma E$. The result is the
$E$-component of the {\em second connection} of
$(M,g,{\cal F})$ \cite{{V2},{V1}}, and we call $\nabla$ the {\em canonical
connection} of $(M,{\cal F},g)$.
If the transversal projective structure
defined by the canonical connection $\nabla$ is projectable, the foliation
$\hat{\cal F}$ of ${\cal Z}({\cal F})$ will have a projectable, almost
complex
structure ${\cal I}_{1}$. But, ${\cal I}_{2}$ will be projectable iff
$\nabla$ is projectable (Proposition 3.4)
hence, ${\cal F}$ is a Riemannian foliation with
the {\em bundle-like metric} $g$ \cite{Mol}.
\section{The integrability conditions}
Now, in the cases where projectability holds
we will discuss
the integrability conditions of the structures ${\cal I}_{1},{\cal I}_{2}$
defined in Section 3 i.e., the conditions which ensure that the local,
almost
complex slice spaces of
$({\cal Z}({\cal F}),\hat{\cal F})$ are, in fact, complex
manifolds.

The integrability
condition is equivalent to the fact that
the local equations of $\hat{\cal F}$ can be put under the form
$z^{\sigma}=const.,
\bar z^{\sigma}=const.$, where $z^{\sigma}$ are complex local
coordinates with complex analytic transition functions. Equivalently, the
distribution ${\cal C}'_{1}$ or ${\cal C}'_{2}$ defined by (3.1$'$),
(3.2$'$) must be {\em Nirenberg integrable}, which (by a theorem of
Nirenberg \cite{Nir}, and since ${\cal C}'_{a}+\bar{\cal C}'_{a}=
T{\cal Z}({\cal F})$ $(a=1,2)$) happens iff ${\cal C}'_{1}$
or ${\cal C}'_{2}$, respectively, is an involutive distribution.

Accordingly,
the structure equations (3.3), (3.4) tell us that the integrability
conditions of ${\cal I}_{1}$ are
$$d\beta=-\omega\wedge\beta-\theta\wedge\bar\beta+T=0\hspace{3mm}(mod.\,
\beta,\theta),\leqno{(5.1)}$$
$$d\theta=-\omega\wedge\theta-\theta\wedge\bar\omega+\Theta=0\hspace{3mm}(mod.\,
\beta,\theta),\leqno{(5.2)}$$
The integrability conditions of
${\cal I}_{2}$ are given by (5.1)
modulo $\beta,\bar\theta$ and (5.2) modulo $\bar\beta,\theta$.
\proclaim 5.1 Theorem. Let $M$ be a foliated manifold with the foliation
${\cal F}$ of codimension $2q$, and assume that $\nabla$ is
a strongly suitable
connection on $N{\cal F}$. Then, ${\cal I}_{2}$ is never integrable, and
${\cal I}_{1}$ is integrable iff
either i) $q\geq3$ and
the torsion and curvature of $\nabla$
are of the following form
$$T_{\nabla}(X,Y)=(\alpha\wedge\pi_{N})(X,Y)=\alpha(X)\pi_{N}(Y)
-\alpha(Y)\pi_{N}(X), \leqno{(5.3)}$$
$$R_{\nabla}(X,Y)(\pi_{N}Z)=\beta(X,Z)(\pi_{N}Y)-\beta(Y,Z)(\pi_{N}X)
+\gamma(X,Y)(\pi_{N}Z),\leqno{(5.4)}$$
or ii) $q=1$, or iii) $q=2$,
the torsion is given by (5.3), and
$$R_{\nabla}(X,Y)\pi_{N}X=\xi(X,Y)\pi_{N}X+\zeta(X,X)\pi_{N}Y,
\leqno{(5.5)}$$
or, equivalently by polarization,
$$R_{\nabla}(X,Y)\pi_{N}Z+R_{\nabla}(Z,Y)\pi_{N}X\leqno{(5.5')}$$
$$=\xi(X,Y)\pi_{N}Z+\xi(Z,Y)\pi_{N}X+2\zeta(X,Z)\pi_{N}Y.$$
In these formulas,
$X,Y,Z$ are vector fields on $M$, $\alpha$ is a $1$-form which vanishes
on $F$, $\gamma$ is a $2$-form which satisfies
$\gamma(X,Y)=\alpha(X,Y)$ with $\alpha$ of (3.6)$_{3}$
if $Y\in F$,
$\beta$ is a $2$-covariant tensor field which vanishes if one of its
arguments is in $F$,
$\xi$ is a two-form such that
$\xi(X,Y)=\alpha(X,Y)$ with $\alpha$ of (3.6)$_{3}$ if $Y\in F$,
and $\zeta$ is a covariant symmetric tensor field which vanishes if one of
its arguments is in $F$. Moreover, after a choice of the
transversal bundle $E$, the values of $\alpha,\beta,\gamma,\xi,\zeta$
are completely determined by traces of $T$ and $R$. \par
\noindent{\bf Proof.} Equation (5.1) never holds
modulo $\beta,\bar\theta$ hence, ${\cal I}_{2}$
is never integrable, and
we only have to discuss the integrability of ${\cal I}_{1}$.

The conclusions could
be justified as follows. By Proposition 3.2, it suffices
to write down the integrability conditions of
the structure ${\cal I}_{1}$
of the local projectable connections $\nabla^{\alpha}$, semisymmetrically
related with $\nabla$. But, on arguments in $N{\cal F}$, these conditions are the
same as after projection on the local slice spaces i.e.,
conditions for twistor spaces on manifolds, where they
exactly are (5.3), (5.4), (5.5) for $\nabla^{\alpha}$
\cite{OR}. Then, it is enough to show that the conditions are invariant by
a semisymmetric transformation of the connection, which is easy.

However, we prefer to give a different, elementary proof.
First, we fix the bundle $E$, and identify
the projections $\pi_{N}$ and $\pi_{E}$. This will produce the
values of the coefficients
$\alpha,\beta,\gamma,\xi,\zeta$ on arguments in $E$, and, then,
the addition of the conditions stated by the theorem
for these coefficients will allow us to extend
$\alpha,\beta,\gamma,\xi,\zeta$ to
arguments in $F$ too.

Equation (5.1)
provides a {\em torsion integrability condition}, and (5.2) provides
a {\em curvature integrability condition}. From (3.7), we see that
(5.1) holds iff $T^{a}_{\bar b\bar c}=0$ for all the frames in
${\cal B}(M,{\cal F})$ (notation of Section 2).
Using again the notation of the proof of Theorem
3.1, we see that the meaning of the previous condition is that
$\forall S^{q}\subseteq E^{c}$ of real index zero and $\forall V,W\in S$,
one must have
$$T_{\nabla}(\bar V,\bar W)\in \bar S.\leqno{(5.6)}$$
As in the proof of Theorem 3.1, if $q\geq2$,
for  any independent vectors $V,W\in E^{c}$ with
$L:=span\{V,W\}$ of real index zero, we may write $L$ as
$L=S_{1}\cap S_{2}$ for some
$S_{1},S_{2}$ as in (5.6). Therefore,
after a conjugation, (5.6) holds iff
$$T_{\nabla}(V,W)=a(V)W+b(W)V.\leqno{(5.7)}$$
Then, because of the skew symmetry, (5.7) becomes (5.3).

Now, from (3.8), and since (3.10) is supposed to hold,
we see that the curvature integrability condition
(5.2) means $R^{a}_{\bar b\bar c\bar d}=0$, and, with the same
notation as above, this is equivalent to
$R_{\nabla}(\bar V,\bar W)(\bar U)\in\bar S$,
whenever the space $L:=span\{V,W,U\}$ has real index zero.
If $q\geq3$, and seeing $L$
as $L=\bar S_{1}\cap\bar S_{2}$,
the integrability condition becomes
$$R_{\nabla}(V,W)U=\gamma(V,W)U+
\beta(V,U)W+\eta(W,U)V.\leqno{(5.8)}$$
Using the skew symmetry in $V,W$, we may write condition (5.8)
under the form (5.4).

For $q=1$, the integrability conditions are satisfied for any strongly
suitable connection.

For $q=2$, (5.3) continues to be
the necessary and sufficient torsion integrability condition (as we saw
above).
Furthermore, (5.4) continues to be a sufficient curvature
integrability condition but, not a necessary one. Indeed, if $q=2$, the
curvature integrability condition
is that for all independent $V,W\in E^{c}$ which span
a space
of real index zero one has
$$R_{\nabla}(V,W)(\lambda V+\mu W)=
\sigma V+\tau W \leqno{(5.9)}$$
$(\lambda,\mu,\sigma,\tau\in{\bf C})$. For (5.9) to hold, it is
enough that
$$R_{\nabla}(V,W)V=\alpha(V,W)V+
\beta(V,V)W.\leqno{(5.10)}$$
(Of course, in (5.9), (5.10) all the coefficients are new.)
Since the condition put on $V,W$ is equivalent with the fact that
the real and imaginary parts
$(Re\,V,Im\,V,Re\,W,Im\,W)$ are ${\bf R}$-independent, we may replace them
by the conjugate vectors as well, and
(5.10) also
provides the necessary and sufficient
curvature integrability condition in terms of real
vectors of $TM$. Namely, we exactly obtain condition (5.5).

Furthermore,
computing with (5.3) the trace of the well defined operator
$i(X)T_{\nabla}:N{\cal F}\rightarrow N{\cal F}$, $X\in\Gamma TM$, we get
$$\alpha(X)=\frac{1}{2q-1}tr\,i(X)T_{\nabla}.\leqno{(5.11)}$$

Similarly, from (5.4) we get
$$tr\,R_{\nabla}(X,Y)=2[alt(\beta)(X,Y)+q\gamma(X,Y)],\leqno{(5.12)}$$
where $alt(\beta)$ is the skew-symmetric part of $\beta$,
valid for all $X,Y\in\Gamma TM$ because of the conditions on $\beta,\gamma$,
and then, with the fixed bundle $E$,
$$Ric_{(\nabla,E)}(X,Y)=
-(2q-1)\beta(X,Y)-\gamma(X,Y), \leqno{(5.13)}$$
valid $\forall X,Y\in\Gamma E$.
By symmetrization, (5.13) gives
$$sym(\beta)(X,Y)=-\frac{1}{2q-1}sym(Ric_{(\nabla,E)})(X,Y),
\hspace{3mm}X,Y\in\Gamma E,\leqno{(5.14)}$$
where $sym$ denotes the symmetric part of a tensor.
Furthermore, the alternation of (5.13) together with (5.12) form an
equation system which yields
$$\gamma(X,Y)=\frac{1}{2q^{2}-q-1}[alt(Ric_{(\nabla,E)})(X,Y)\leqno{(5.15')}$$
$$+\frac{2q-1}{2}
tr\,R_{\nabla}(X,Y)],\hspace{3mm}X,Y\in\Gamma E,$$
$$alt\,\beta(X,Y)=-\frac{1}{2q^{2}-q-1}[q\,alt\,Ric_{(\nabla,E)}(X,Y)
+\frac{1}{2}tr\,R_{\nabla}(X,Y)].\leqno{(5.15'')}$$
Finally,
the value of $\beta$ is provided by inserting the values (5.15$'$),
(5.15$''$) into
$\beta=sym(\beta)+alt(\beta)$.

For the case $q=2$, we compute traces in the $E$-version of
(5.5$'$) and, first, $\forall X,Y\in\Gamma E$, we get
$$tr\,R_{\nabla}(X,Y)+Ric_{(\nabla,E)}(Y,X)=
5\xi(X,Y)+2\zeta(X,Y).\leqno{(5.16)}$$
Then, we alternate and symmetrize, and get
$$\xi(X,Y)=\frac{1}{5}[tr\,R_{\nabla}(X,Y)-alt(Ric_{(\nabla,E)})(X,Y))],
\leqno{(5.17)}$$
$$\zeta=\frac{1}{2}sym(Ric_{(\nabla,E)})(X,Y),
\hspace{3mm}X,Y\in\Gamma E.\leqno{(5.18)}$$
Q.e.d.

Now, let us discuss the integrability conditions
(5.1), (5.2)of the projectable
structure ${\cal I}_{1}$ defined by an adapted torsionless connection
$\nabla$ which is projectively related to a given transversal projective
structure of the foliation ${\cal F}$ (Theorem 4.3).

Condition (5.1) holds since $T=0$. In
(5.2), $\Theta$ has one more term than in the case of Theorem 5.1 but,
this is a term in $\beta^{c}\wedge\gamma^{u}$ (see (3.8)).
Therefore, the integrability condition (5.2) again means
$R^{a}_{\bar b\bar c\bar d}=0$, as in Theorem 5.1.
Hence, the only integrability condition for
${\cal I}_{1}$ is again (5.4), if $q\geq3$, and (5.5), if $q=2$.

Furthermore, since there is no torsion,
we have the transversal Bianchi identity
$$\sum_{Cycl(X,Y,Z)}R_{\nabla}(X,Y)Z=0,\hspace{5mm}X,Y,Z\in\Gamma E,
\leqno{(5.19)}$$
and, if we insert the value of $R_{\nabla}$ as given by the $E$-version of
(5.4) in (5.19), and take the trace, we obtain
$\gamma(X,Y)=\beta(X,Y)-\beta(Y,X)$, $\forall X,Y\in\Gamma E$.

The final result can be expressed as follows
\proclaim 5.2 Theorem. Let ${\cal F}$ be a transversally projective
foliation of codimension $2q$, $q\geq3$. Then its family of global, adapted,
projectively equivalent, torsionless connections of $N{\cal F}$ defines a
family of projectable structures ${\cal I}_{1}$ which are simultaneously
integrable or not. Integrability occurs iff
$$R_{\nabla}(X,Y)\pi_{N}Z=\beta(X,Z)\pi_{N}Y-\beta(Y,Z)
\pi_{N}X\leqno{(5.20)}$$
$$+(\beta(X,Y)-\beta(Y,X))\pi_{N}Z,\hspace{3mm} X,Y,Z\in\Gamma TM, $$
where necessarily,
$\beta(X,Y)=0$ if $Y\in\Gamma F$, $\beta(Y,X)=-\lambda(X,Y)$,
with $\lambda$ of (3.28), if $Y\in\Gamma F$, and for any fixed $E$, $\forall
X,Y\in\Gamma E$,
$$\beta(X,Y)=\frac{1}{2q-1}sym(Ric_{(\nabla,E)})(X,Y)+
\frac{1}{4q+2}tr\,R_{\nabla}(X,Y). \leqno{(5.21)}$$
The meaning of condition (5.20) is that the transversal projective structure
of ${\cal F}$ is projectively flat. Under the same hypotheses, but for
$q=2$, ${\cal I}_{1}$ is integrable iff
$$R_{\nabla}(X,Y)\pi_{N}X=\xi(X,Y)\pi_{N}X\hspace{5mm}X,Y\in\Gamma
TM,\leqno{(5.22)}$$
where $\xi$ is a $2$-form which satisfies $\xi(X,Y)=2\lambda(X,Y)$,
with $\lambda$ of (3.28), if $Y\in\Gamma F$, and if $E$ is fixed,
$$\xi(X,Y)=\frac{1}{5}[tr\,R_{\nabla}(X,Y)-Ric_{(\nabla,E)}(X,Y)],
\hspace{3mm}X,Y\in\Gamma E.\leqno{(5.23)}$$
For $q=1$, integrability always holds.\par
\noindent{\bf Proof.} Condition (5.20) was already justified, and (5.21)
follows by taking the necessary traces in (5.20). Condition (5.22) has a
similar justification namely, the use of (5.5$'$) and (5.19) yields
$\zeta=0$ hence, we have proven (5.22) for arguments in $\Gamma E$.
Then, (5.23) follows from (5.17), (5.18), and $\zeta=0$.

Finally, the relation which a projective transformation (4.6) of a connection
implies on the curvature tensors can be easily deduced, and it is well known
\cite{Kob}. This relation preserves the form (5.20) of the curvature
(while changing $\beta$ of course), and (5.20) holds iff the connection is
projectively transformable into a flat connection.
This explains all the other assertions of the theorem. Q.e.d.

We may look for examples by first defining the notion of
a {\em transversally homographic foliation} as follows. In  ${\bf
R}^{s}$, there exists a pseudogroup of {\em homographies},
i.e., projective transformations written in non homogeneous
coordinates
$$\tilde x^{i}=\frac{a^{i}_{h}x^{h}+b^{i}}{c^{i}_{h}x^{h}+d^{i}},
\leqno{(5.24)}$$
where $(i,h=1,...,s)$ and
the coefficients are real numbers; homographies are defined on
open sets of ${\bf R}^{s}$.
A foliation ${\cal F}$ of codimension $s$ will be
called {\em transversally homographic} if its Haefliger cocycle \cite{Mol}
may be taken in the pseudogroup (5.24).

Since homographies send a line segment onto a line segment, the system
of local adapted connections of $N{\cal F}$ defined by the flat connection
of ${\bf R}^{s}$ consists of projectively related connections. Hence, every
transversally homographic foliation is transversally projective flat.
Clearly, the converse also holds. Thus, the two notions coincide.

A concrete example is the following one given in \cite{NS}. Let $G$ be the
connected component of the identity
in the general projective group in dimension $s+1$, $H$
be the affine subgroup of $G$, which fixes a given hyperplane, and $K$ the
subgroup of $H$ which is isomorphic to $SO(s)$. Then $B:=G/H$ is
diffeomorphic to the real projective space ${\bf RP}^{s}$, and
$V:=G/K$ is a bundle over $B$. The fibers of this bundle define a
$G$-invariant, transversally homographic foliation ${\cal S}$
of codimension $s$. Then, ${\cal
S}$ induces transversally homographic foliations ${\cal F}$ on all the
manifolds $M=D\backslash V$ where $D$ is a discrete subgroup of $G$.

In particular, if $s=2q$, all the foliations ${\cal F}$ obtained in this way
have a twistor space ${\cal Z}({\cal F})$ with a transversally holomorphic
lifted foliation $\hat{\cal F}$.

Other examples could be provided by foliations ${\cal F}$ on
Riemannian manifolds $(M,g)$, such that the canonical connection of the
normal bundle (see the end of Section 4) is projectively flat but, we do
not have yet concrete examples of such foliations.
\section{The case of a Riemannian manifold}
In this section we discuss the transversal twistor space ${\cal Z}({\cal
F})$ of a foliation ${\cal F}$ defined on a Riemannian manifold $M$ with
the metric $g$.

Before doing this
we notice that a general metric $g$ may not allow for a
twistorial construction,
in the sense of Section 2, on the bundle ${\cal Z}({\cal F},g)$ of
pointwise, $g/_{E}$-orthogonal, complex structures of $E$, because
${\cal Z}({\cal F},g)$ may not be a foliated bundle, and then the foliation
$\hat{\cal F}$ does not exist. If ${\cal F}$ is a Riemannian foliation with
the {\em bundle-like metric} $g$, a twistor bundle ${\cal Z}({\cal F},g)$
with a foliation $\hat{\cal F}$
exists. For this bundle we have
$${\cal Z}({\cal F},g)
=\cup{\cal Z}({\cal F}/_{U_{\alpha}},g_{\alpha}=g/_{U_\alpha}),\leqno{(6.1)}$$
where $\{U_{\alpha}\}$ is a covering of $M$ by ${\cal F}$-adapted coordinate
neighborhoods, and the transversal almost complex structures
of $\hat{\cal F}$ defined by the Levi-Civita connections of
$g_{E}:=g/_{E}$ \cite{Mol}
are the lifts of those of the Riemannian twistor spaces of the slice
spaces $(U_{\alpha}/{\cal F}\cap U_{\alpha},g_{E}/_{U_{\alpha}})$. Hence,
all the results of a local character (e.g., integrability)
will be just the same as the classical ones e.g.,
\cite{{AHS},{OR},{JR}}. Moreover,
if ${\cal F}$ is a {\em conformal foliation} i.e., only the conformal
structure defined by $g_{E}$ is projectable \cite{{NS},{V3}}, a
space (6.1)
where $g_{\alpha}$ are foliated and conformal to $g_{E}/_{U_\alpha}$,
exists. This is the space of {\em conformal twistors},
and it has a well
defined almost complex structure. Namely, the lift of the conformal
invariant, almost complex structure of the local slice spaces
$(U_{\alpha}/{\cal F}\cap U{_\alpha},g_{\alpha})$
\cite{JR}.

Now, coming back to the space ${\cal Z}({\cal F})$, the interesting point
to be made is that, using an adapted connection $\varpi$ of ${\cal F}$,
$g$ can be lifted to a pseudo-Riemannian metric $\hat g$ on
${\cal Z}({\cal F})$ as follows. On the principal bundle ${\cal B}(M,{\cal
F})$ of Section 2, we can define the following matrices of globally defined
functions
$$g_{1}=\,^{t}C\cdot C,\:g_{2}=\,^{t}b\cdot \bar b=\,^{t}\bar g_{2}\:\:
(b,\bar b, C)\in{\cal B}(M,{\cal F}), \leqno{(6.2)}$$
where the dot is the $g$-scalar product at $\pi(b,
\bar b,C)\in M$, and $t$ denotes
transposition of matrices. Then, it is an immediate consequence
of formulas (2.17), (2.18) that
$$\hat
g:=\,^{t}\gamma\otimes(g_{1}\gamma)+\,^{t}\beta\otimes(g_{2}\bar\beta)
+\,^{t}\bar\beta\otimes(\,^{t}g_{2}\beta)
-tr(\theta\otimes\bar\theta+\bar\theta\otimes\theta),\leqno{(6.3)}$$
where the notation is that of Proposition 2.1, pulls back to a well defined
pseudo-Riemannian metric of ${\cal Z}({\cal F})$ of $(\pm)$-signs
$(p+2q+q^{2},q^{2})$.
\proclaim 6.1 Remark. The metric $\hat g$ also exists if $g$ is a
pseudo-Riemannian metric which is nondegenerate on $F$ but, the signature
will be different. \par
It follows from (6.3) that $\hat g/_{\hat E}$ is compatible with the almost
complex structures ${\cal I}_{1},{\cal I}_{2}$ defined by the same adapted
connection $\varpi$. Accordingly, there exists a real $2$-form of
$\hat{\cal F}$-type $(2,0)$ on ${\cal Z}({\cal F})$, given by
$$\Phi=\sqrt{-1}\{\,^{t}\beta\wedge(g_{2}\bar\beta)
-tr(\theta\wedge\bar\theta)\},\leqno{(6.4)}$$
which is nondegenerate on $\hat E$.
\proclaim 6.2 Remark. For the construction of a $2$-form $\Phi$ as in
(6.4), it suffices to start with a pseudo-Riemannian transversal metric of
${\cal F}$ on $M$.\par
\proclaim 6.3 Proposition. The $2$-form $\Phi$ is $\hat{\cal
F}$-projectable iff $g$ is an ${\cal F}$-bundle-like metric on $M$, and
$\varpi$ is a strongly twistor-suitable connection. In particular, if
$\varpi$ is the canonical connection of ${\cal F}$, $\Phi$ is projectable
iff $g$ is a bundle-like metric of ${\cal F}$.\par
\noindent{\bf Proof.} From $\hat{\cal F}$-type considerations and formula
(3.9), it follows that $d''\Phi=0$ iff each term of the right hand side of
(6.4) is $d''$-closed. For the second term this happens iff (3.10) and
(3.11) hold i.e., $\varpi$ is a strongly suitable connection. For the first
term in (6.4), if we put $b_{i}=\xi^{a}_{i}X_{a}$, where $\{X_{a}\}$ is a
local foliated basis of $E$,
then $(\xi_{i}^{a})$ are $\hat{\cal F}$-local
transversal coordinates on ${\cal Z}({\cal F})$, and $d''$-closedness holds
iff $g(X_{a},X_{b})$ are foliated functions on $M$. Q.e.d.
\proclaim 6.4 Corollary. If the conditions of Proposition 6.3 hold, $\hat
g$ is an $\hat{\cal F}$-bundle-like metric.\par
We may also notice that the forms defined on ${\cal B}({\cal F})$ by
$$\hat h:=\,^{t}\beta\otimes(g_{3}\beta)+\,^{t}\bar\beta\otimes(\bar g_{3}
\bar\beta), \leqno{(6.5)}$$
$$\hat\Psi:=
\sqrt{-1}[\,^{t}\beta\wedge(g_{3}\beta)+\,^{t}\bar\beta\wedge(\bar g_{3}
\bar\beta)], \leqno{(6.6)}$$
where $g_{3}:=\,^{t}b\cdot b$, descend to well defined global forms on
${\cal Z}({\cal F})$. (In fact, each term of (6.5) and (6.6) has this
property.) Then, $\hat g+\hat h$ is again a pseudo-Riemannian metric on
${\cal Z}({\cal F})$, whose ${\cal H}$-term is the lift of $g$, and
$\Phi+\Psi$ is an $\hat{\cal F}$-transversal almost symplectic
structure which is projectable iff the conditions of Proposition 6.3 hold.
 \vspace*{1cm}
{\small Department of Mathematics, \\}
{\small University of Haifa, Israel\\}
{\small E-mail: vaisman@math.haifa.ac.il\\}

\begin{thebibliography} {99}
\bibitem{AHS} Atiyah M. F., Hitchin N. J. and Singer I. M., Self-duality in
four-dimensional Riemannian geometry. Proc. Royal Soc. London, A 362 (1978),
425-461.
\bibitem{Bt} R. Bott, Lectures on characteristic classes of foliations.
Lecture Notes in Math. 279, Springer Verlag, Berlin, 1972, 1-94.
\bibitem{JR} Jensen G. R. and Rigoli M., Twistor and Gauss lifts of
surfaces in four-manifolds. Contemporary Math., American Math. Soc., 132
(1989), 197-232.
\bibitem{KT} F. Kamber and Ph. Tondeur, Foliated bundles and
characteristic classes. Lecture Notes in Math. 493, Springer Verlag,
Berlin, 1975.
\bibitem{Kob} S. Kobayashi, Transformation groups in differential
geometry, 2$^{\rm nd}$ ed. Springer Verlag, Berlin, 1995.
\bibitem{KN} Kobayashi S. and Nomizu K., Foundations of Differential
Geometry I, II. Intersci. Publ., New York, 1963, 1969.
\bibitem{Mol2} P. Molino, Propri\'et\'es cohomologiques et propri\'et\'es
topologiques des feuilletages \`a connexion transverse projetable. Topology,
12 (1973), 317-325.
\bibitem{Mol} P. Molino, Riemannian Foliations. Progress in Math. Series,
73, Birkh\"auser, Boston, 1988.
\bibitem{Nir} L. Nirenberg, A complex Frobenious theorem. In:
Sem. on Analytic Functions, Vol. I, Inst. for Adv. Study, Princeton, 1957,
172-179.
\bibitem{NS} S. Nishikawa and H. Sato, On characteristic classes of
Riemannian, conformal and projective foliations. J. Math. Soc. Japan,
28 (1976), 223-241.
\bibitem{OR} O'Brian N. R. and Rawnsley J. H., Twistor Spaces. Annals of
Global Analysis and Geometry, 3 (1985), 29-58.
\bibitem{Pn} Penrose R., The twistor programme. Reports on Math. Physics,
12 (1977), 65-76.
\bibitem{V2} I. Vaisman, Vari\'et\'es Riemanniennes
Feuiellet\'ees. Czechosl. Math. J., 21 (1971), p. 46-75.
\bibitem{V1} I. Vaisman, Cohomology and Differential Forms. M. Dekker,
Inc., New York, 1973.
\bibitem{V3} I. Vaisman, Conformal Foliations. Kodai Math. J., 2(1979),
26-37.
\bibitem{V4} I. Vaisman, A Note on Projective Foliations. Publ. Univ.
Autonoma Barcelona, 27 (1983), 109-128.
\bibitem{V5} I. Vaisman, On the twistor space of almost Hermitian manifolds.
Annal. Global Anal. Geom., 16 (1998), 335-356.
\end{thebibliography}
\end{document}